\documentclass[12pt]{article}

\usepackage{latexsym,amssymb,graphics}
\usepackage{color}
\usepackage{epsfig}

\def\tr{^{\rm T}}
\def\Real{\mathbb R}

\def\Integer{\mathbb N}

\def\dst{\displaystyle}

\def\qmx#1{\left(\matrix{#1}\right)}

\def\beeq#1{\begin{equation}{#1}\end{equation}}

\def\ba{\begin{array}}
\def\ea{\end{array}}
\def\eqa{\begin{eqnarray}}
\def\eqe{\end{eqnarray}}

\newtheorem{proposition}{Proposition}

\newtheorem{lemma}{Lemma}

\newenvironment{proof}{\medskip\noindent{\it Proof. }}{ \medskip}
\newenvironment{remark}{\medskip\noindent{\it Remark. }}{
\medskip}

\begin{document}
\title{Remote Tracking via Encoded Information\\ for Nonlinear Systems\thanks{This work was partially
supported by NSF under grant ECS-0314004, by ONR under grant
N00014-03-1-0314, and by MIUR. Corresponding Author: Dr. Lorenzo
Marconi, email: lmarconi@deis.unibo.it, tel. 0039 051 2093788,
fax. 0039 051 2093073}}

\author{A. Isidori $^{\dag \circ\ddag }$ ,  L. Marconi $^{\circ}$,
C. De Persis $^{\dag}$}
\date{\today}

\maketitle
\begin{center}

\small $^{\dag}$Dipartimento di Informatica e Sistemistica,
 Universit\`{a} di Roma ``La Sapienza'', \\00184 Rome, ITALY.

$^{\circ}$ C.A.SY. -- Dipartimento di Elettronica, Informatica e
Sistemistica, University of Bologna, \\40136 Bologna, ITALY.

$^{\ddag}$Department of Electrical and Systems Engineering,
Washington University, \\St. Louis, MO 63130.
\smallskip
\normalsize
\end{center}
\maketitle
\begin{abstract}
The problem addressed in this paper is to control a plant so as to
have its output tracking (a family of) reference commands
generated at a {\em remote location} and transmitted through a
communication channel of finite capacity. The uncertainty
due to the presence of the communication
channel is counteracted by a suitable choice of the parameters of the regulator.
\end{abstract}


\section{Introduction}

In distributed control systems, sensors, actuators and control unit may be
placed at locations which are geographically separated. Information among these
devices must then be exchanged through a finite
bandwidth channel.\\
The problem addressed in this paper is to control a plant so as to
have its output tracking (a family of) reference commands
generated at a {\em remote location} and transmitted through a
communication channel of finite capacity. What renders the problem
in question different from a conventional tracking problem is that
the {\em tracking error}, that is the difference between the
command input and the controlled output, is not available as a
{\em physical entity}, as it is defined as difference between two
quantities residing at different (and possibly distant) physical
locations. Therefore the tracking error as such cannot be used to
drive a feedback controller, as it is the case in a standard
tracking problem.\\
The actual tracking error not being available, it is natural to
approach the problem by reconstructing the tracking error starting from the information
transmitted through the communication
channel. For the reconstruction to successfully take place, the
information must be suitably encoded. A possibility  is
to make use of a sufficiently large number of bits, so as to render
the magnitude of the difference between the true reference signal and the
reconstructed one negligible. However, in the framework of distributed control
systems, the constraint on the available bandwidth is usually tight and adopting
encoding schemes which require a large number of bits may not be
practically feasible. The approach pursued in this paper is rather
to counteract the uncertainty due to the presence of the communication
channel by a suitable choice of the parameters of the regulator.\\
In Section 2, the formulation of the problem is made more precise,
whereas the procedure for encoding the reference command is described in
Section 3. Section 4 introduces the regulator which guarantees the
achievement of the control goal using the reconstructed tracking
error. The main results of the paper are stated and proved in Section 5.
The proof consists of two steps. First, boundedness of the
closed-loop trajectories are shown (Section 5.1), and then asymptotic
convergence to zero of the tracking error is concluded (Section
5.2). The results are illustrated by an example in Section 6.

\section{Problem statement}

Generally speaking, the problem in question can be defined in the
following terms. Consider a single-input single-output nonlinear
system modeled by equations of the form \beeq{\label{sysgeneral}
\ba{rcl} \dot x &=& f(x)+g(x)u\\
y &=& h(x)\ea} and suppose its output $y$ is required to
asymptotically track the output $y_{\rm des}$ of a remotely
located exosystem
\beeq{\label{exogeneral}
\ba{rcl} \dot w &=& s(w) \qquad w \in \Real^r\\
y_{\rm des} &=& y_{\rm r}(w)\,.\ea}

The problem is to design a control law of the form
\beeq{\label{contrgeneral}
\ba{rcl} \dot \xi &=& \varphi(\xi,y,w_{\rm q})\\
u &=& \theta(\xi,y,w_{\rm q})\ea} in which $w_{\rm q}$ represents
a {\em sampled and quantized} version of the remote exogenous
input $w$, so as to have the tracking error
\beeq{\label{trackerror} e(t) = y(t) - y_{\rm r}(w(t))}
asymptotically converging to zero as time tends to $\infty$. Note
that the controller in question does not have access to $e$, which
is not physically available, but only to the controlled output and
to a sampled and quantized version of the remotely generated
command.

We will show in what follows how the theory of output tracking can
be enhanced so as to address this interesting design problem. In
particular, we will show how, by incorporating in the controller
two (appropriate) internal models of the exogenous signals, the
desired control goal can be achieved. One internal model is meant
to asymptotically reproduce, at the location of the controlled
plant, the behavior of the remote command input. The other
internal model, as in any tracking scheme, is meant to generate
the ``feed-forward" input which keeps the tracking error
identically at zero.

We begin by describing, in the following section, the role of the
first internal model.

 \section{The encoder-decoder pair} \label{secEncDec}
 In order to overcome the limitation due to the finite capacity of
 the communication channel, the control structure proposed here
 has a decentralized structure consisting of two separate units:
 one unit, co-located with the command generator,
  consists of an {\em encoder} which extracts from the
  reference signal the data which are transmitted
 through the communication channel; the other unit, co-located
 with the controlled plant, consists of a  {\em decoder}
which processes the encoded received information and of a {\em
regulator} which generates appropriate control input.

 The problem at issue will be solved under a number of assumptions
 most of which are inherited by the literature of output regulation
 and/or control under quantization. The first assumption, which is
 a customary condition in the problem of output regulation, is
 formulated as follows.

\medskip\noindent{\bf (A0)}  The vector field $s(\cdot)$ in
(\ref{exogeneral}) is locally Lipschitz and the initial conditions
for (\ref{exogeneral}) are taken in a fixed compact invariant set
$W_0$.$\;\;\triangleleft$

 \medskip The next assumption is, on the contrary, newer and
 motivated by the specific problem addressed in this paper.
 In order to formulate rigorously the assumption in question,
 we need to introduce some notation.
 In particular let $|x|_S$ denote the distance  at a point $x \in \Real^n$
 from a compact subset $S \subset \Real^n$, i.e. the number
 \[
 |x|_S := \max_{y \in S} |x - y|
 \]
 and let
 \beeq{ \label{L0def}
 L_0 = \max_{i \in [1, \ldots, r] \atop (x,y) \in W_0 \times W_0} |x_i -
 y_i|\,.
 }

 Furthermore, having denoted by $N_b$ the number of bits
 characterizing the communication channel constraint, let $N$ be
 the largest positive integer such that
 \beeq{\label{nbNrel}
 N_b \geq  r \, \lceil \log_2N\rceil
 }
 where $\lceil \upsilon \rceil$, $\upsilon \in \Real$, denotes the lowest integer such
 that $\lceil \upsilon \rceil \geq \upsilon$.

 With this notation in mind, the second assumption can be
 precisely formulated as follows.

\medskip\noindent{\bf (A1)} There exists a compact set $W \supset W_0$
 which is invariant for $\dot w = s(w)$ and such that
 \[
  \bar w \not\in W  \qquad \Rightarrow \qquad |\bar w|_{W_0}
  {>} \sqrt{r} {L_0 \over 2 N}\,. \;\;\;
 \triangleleft
 \]

 \medskip $W$ being compact and $s(\cdot)$ being locally Lipschitz,
 it is readily seen that there exists a non decreasing and bounded function
 $M(\cdot): \Real_{\geq 0} \rightarrow \Real_{>0}$, with $M(0) = 1$, such that for all
 $w_{10} \in W$ and $w_{20} \in W$ and for all $t \geq 0$
 \beeq{\label{Mdef}
 |w_1(t) - w_2(t)| \leq M(t) |w_{10} - w_{20}|
 }
 where $w_1(t)$ and $w_2(t)$ denote the solutions of (\ref{exogeneral})
 at time $t$ passing through $w_{10}$ and, respectively, $w_{20}$ at time $t=0$.

 This function, the sampling interval $T$, the number $L_0$ defined in (\ref{L0def})
 and the number $N$ fulfilling (\ref{nbNrel}), determine the parameters of the encoder-decoder
pair, which are defined as follows (see \cite{TaMi}, \cite{LiHe}, \cite{DePI}
for more details).
\medskip

 \noindent {\em Encoder dynamics}. The encoder dynamics consist
 of a copy of the exosystem dynamics, whose state is updated at
 each sampling time $k T$, $k \in  \Integer$, and determines
 (depending on the actual state of the exosystem) the centroid of the
 quantization region, and of an
 additional discrete-time dynamics which determines
  the size of the quantization region. Specifically, the encoder
  is characterized by
 \[
 \ba{rcll}
 \dot w_{\rm e} &=& s(w_{\rm e})&  \qquad
 w_{\rm e}(kT) = w_{\rm e}(kT^-) + w_{\rm q}(k)\dst {L(k) \over N} \qquad w_{\rm
 e}(0^-) \in W_0
 \\[2mm]
 L(k+1) &=& \dst \sqrt{r} \, {M(T) \over  N} L(k) & \qquad L(0) =
 L_0
 \ea
 \]
 in which $w_{\rm q}$ represents the encoded information
 given by, for $i=1,\ldots,r$,
 \[
 w_{{\rm q},i}(k) = \mbox{sgn}(w_i(kT)- w_{e,i}(kT^-)) \cdot \left
 \{ \ba{ll}
 \dst \left \lceil \dst {N|w_i(kT) - w_{e,i}(kT^-)| \over L(k)} \right \rceil - \dst {1 \over 2} &
 \hspace{5mm} N \;
 \mbox{even}\\[4mm]
 \dst \left \lceil \dst {N|w_i(kT) - w_{e,i}(kT^-)| \over L(k)} - \dst
 {1 \over 2} \right \rceil & \hspace{5mm} N \;
 \mbox{odd}\,.
 \ea \right .
 \]
 At each sampling time $kT$, the vector $w_{\rm q}(k)$ is transmitted to the controlled plant
 through the communication channel and then used to update the state of the decoder
 unit as described in the following. To this regard note that
 each component of the vector $w_{\rm q}(k)$ can be described by $\lceil \log_2N\rceil$
 bits and thus the communication channel constraint  is fulfilled.
 \medskip

 \noindent {\em Decoder dynamics}
 The decoder dynamics is a replica of the encoder dynamics and it
 is given by
  \beeq{\label{decdyn}
 \ba{rcll}
 \dot w_{\rm d} &=& s(w_{\rm d})&  \qquad w_{\rm d}(kT) = w_{\rm d}(kT^-) + w_{\rm q}(k)
 \dst {L(k) \over N} \qquad w_{\rm d}(0^-) = w_{\rm e}(0^-)\\[2mm]
 L(k + 1) &=& \dst \sqrt{r} \, {M(T) \over N} L(k) & \qquad
 L(0) = L_0
 \ea
 }
 \medskip

 If, ideally, the communication
 channel does not introduce delays, it turns out that
 $w_{\rm d}(t) \equiv w_{\rm e}(t)$ for all $t \geq 0$.
 Furthermore, it can be proved that the set $W$ characterized in
 Assumption (A1) is invariant for the encoder (decoder) dynamics
 and that the asymptotic behavior of
 $w_{\rm e}(t)$ ($w_{\rm d}(t)$) converges uniformly to the true
 exosystem state $w(t)$, provided that $T$ is properly
 chosen with respect to the number $N$ and the function $M(\cdot)$. This is formalized
 in the next proposition (see \cite{LiHe},
\cite{DePI} for details).

 \begin{proposition}\label{PR1}
 Suppose Assumptions (A0)-(A1) hold and that the sampling time $T$ and the number $N$
 satisfy
 \beeq{ \label{TNconst}
 N  > \sqrt{r} \, M(T) \,.
 }
 Then:\\[1mm]

 \noindent (i) for any $w_{\rm d}(0^-) \in
 W_0$ and $w(0) \in W_0$, $w_{\rm d}(t) \in W$ for all $t \geq 0$;\\[2mm]

 \noindent (ii) for any $w_{\rm d}(0^-) \in W_0$ and $w(0) \in W_0$,
 \[
\lim_{t \rightarrow \infty} |w(t) - w_{\rm d}(t)| = 0
 \]
 with uniform convergence rate, namely for every $\epsilon>0$ there exists $T^\ast>0$
 such that for all initial states $w_{\rm d}(0^-) \in W_0$, $w(0) \in W_0$,  and for all $ t \geq T^\ast$,
 $|w(t) - w_{\rm d}(t)| \leq \epsilon$.
 \end{proposition}
 \begin{proof}
 As $W$ is an invariant set for $\dot w = s(w)$, the proof of the
 first item reduces to show that, for all $k \geq 0$, if $w_{\rm d}(kT^-) \in W$ then
 necessarily $w_{\rm d}(kT) \in W$.
 For, note that this is true for $k=0$. As a matter of fact,
 since $w_{\rm d}(0^-) \in W_0 \subset W$ and by bearing in mind the definition
 of $w_{\rm q}$, it
 turns out that $|w_{\rm d}(0) - w(0)| \leq \sqrt{r} {L_0 / 2 N}$
 which implies, by definition of $W$ in Assumption (A1), that
 $w_{\rm d}(0) \in W$. For a generic $k>0$ note that, again by
 definition of $w_{\rm q}$, it turns out that
 $|w_{\rm d}(kT) - w(kT)| \leq \sqrt{r} {L(k) / 2 N}$. But, by the second
 of (\ref{decdyn}) and by condition (\ref{TNconst}), $L(k) < L(k-1) \leq L_0$
 yielding $|w_{\rm d}(kT) - w(kT)| \leq \sqrt{r} {L_0 / 2 N}$
 which implies $w_{\rm d}(kT) \in W$. This completes the proof of
 the first item. The second item has  been proved in \cite{LiHe}, \cite{DePI}. $\triangleleft$
 \end{proof}

 \begin{remark}
 By composing  (\ref{nbNrel}) with (\ref{TNconst}) it is easy to realize that
 the number of bits $N_b$ and the sampling interval $T$ are required to satisfy
 the constraint
 \beeq{\label{Nbconst}
 N_b \geq r \left \lceil \log_2 \left ( \sqrt{r} \, M(T)\right ) \right
 \rceil
 }
 in order to have the encoder-decoder trajectories asymptotically converging to the
 exosystem trajectories. Since the function $M(\cdot)$ depends on the exosystem
 dynamics and on the set $W_0$ of initial conditions for (\ref{exogeneral}),
 equation (\ref{Nbconst}) can be interpreted as a relation  between
 the bit-rate of the communication channel and the exosystem dynamics which must
 be satisfied in order to remotely reconstruct the reference signal.
 \end{remark}

 \section{The regulator}

 \subsection{Basic hypotheses}

As in most of the literature on regulation of nonlinear system, we
assume in what follows that the controlled plant has well defined
relative degree and normal form. If this is the case and if the
initial conditions of the plant are allowed to vary on a fixed
(though arbitrarily large) compact set, there is no loss of
generality in considering the case in which the controlled plant
has relative degree 1 (see for instance \cite{IFAC05}). We
henceforth suppose that system (\ref{sysgeneral}) is expressed in
the form \beeq{\label{sys1}
 \ba{rcll}
 \dot z &=& f(z,y,\mu)&\qquad z \in \Real^n\\
 \dot y &=& q(z,y,\mu) + u & \qquad y \in \Real
 \ea
 }
in which $\mu$ is a vector of uncertain parameters ranging in a
known compact set $P$.  Initial conditions $(z(0),y(0))$ of
(\ref{sys1}) are allowed to range on a fixed (but otherwise
arbitrary) compact set $Z
\times Y \subset \Real^n \times \Real$.

It is well known that, if the regulation goal is achieved, in
steady-state (i.e when the tracking error $e(t)$ is identically
zero) the controller must necessarily provide an input of the form
\beeq{\label{qzero} u_{\rm ss} =  L_s y_{\rm r}(w)- q(z,y_{\rm
r}(w),\mu) } (where $L_s y_{\rm r}(\cdot)$ stands for the
derivative of $y_{\rm r}(\cdot)$ along the vector field
$s(\cdot)$) in which $w$ and $z$ obey \beeq{\label{fzero}
\ba{rcl}\dot \mu&=& 0\\
\dot w &=& s(w)\\
\dot z &=& f(z, y_{\rm r}(w),\mu)\,.\ea}

As in \cite{BI03}, we assume in what follows that system
(\ref{fzero}) has a compact attractor, which is also locally
exponentially stable. To express this assumption in a concise
form, it is convenient to group the components $\mu,w,z$ of the
state vector of (\ref{fzero}) into a single vector ${\bf z}={\rm
col}(\mu,w,z)$ and rewrite the latter as
\[
\dot {\bf z} = {\bf f}_0({\bf z})\,.
\]
Consistently, the map (\ref{qzero}) is rewritten as
\[
u_{\rm ss} = {\bf q}_0({\bf z})\,,\] and it is set ${\bf
Z}=P\times W\times Z$. The assumption in question is the following
one\,\footnote{Recall that, if the positive orbit of a compact set
$X$ of initial conditions of a system
\beeq{\label{uno}\dot x=f(x)} is bounded, the $\omega$-limit set
of $X$ under the flow of (\ref{uno}) -- denoted $\omega(X)$ -- is
a nonempty, compact invariant set which attracts $X$ uniformly. If
$\omega(X)$ is in the interior of $X$, then $X$ is also stable in
the sense of Lyapunov.}

\medskip\noindent {\bf (A2)}
there exists a compact subset $\mathcal Z$ of $P\times W \times
\Real^n$ which contains the positive orbit of the set ${\bf Z}$
under the flow of (\ref{fzero}) and $\omega({\bf Z})$ is a
differential submanifold (with boundary) of $P\times W
\times
\Real^n$. Moreover there exists a number $d_1>0$ such that
 \[
{\bf z} \in P\times W \times \Real^n\,, \quad |{\bf
z}|_{\omega({\bf Z})} \leq d_1 \qquad \Rightarrow \qquad {\bf z}
\in {\bf Z}\,.
\]
Finally,  there exist $m \geq 1$, $a>0$ and $d_2 \leq d_1$ such that
\[
{\bf z}_0 \in P\times W \times \Real^n\,, \quad |{\bf
z}_0|_{\omega({\bf Z})} \leq d_2 \qquad \Rightarrow \qquad
 |{\bf z}(t)|_{\omega({\bf Z})} \leq m e^{- a t}
 |{\bf z}_0|_{\omega({\bf Z})}\,,
\]
in which ${\bf z}(t)$ denotes the solution of (\ref{fzero})
passing through ${\bf z}_0$ at time $t=0$. $\triangleleft$

\medskip
In what follows, the set $\omega({\bf Z})$ will be simply denoted
as ${\mathcal A}_0$. The final assumption is an assumption that
allows us to construct an {\em internal model} of all inputs of
the form $u_{\rm ss}(t) = {\bf q}_0({\bf z}(t))$, with ${\bf
z}(t)$ solution of (\ref{fzero}) with initial condition in
${\mathcal A}_0$. This assumption, which can be referred to as
assumption of {\em immersion} into a {\em nonlinear uniformly
observable} system, is the following one.

\medskip\noindent
{\bf (A3)} There exists an integer $d>0$ and a locally Lipschitz
map
 $ \varphi:\Real^d \rightarrow \Real$ such that,  for all ${\bf z} \in
{\mathcal A}_0$, the solution ${\bf z}(t)$ of (\ref{fzero})
passing through ${\bf z}_0$ at $t=0$  is such that the function
$u(t)={\bf q}_0({\bf z}(t))$ satisfies
 \beeq{ \label{nlimc}
  u^{(d)}(t) +
 \varphi(u^{(d-1)}(t), \ldots, u^{(1)}(t),u(t))=0\,.
 \;\;\triangleleft
 }

\subsection{The design of the regulator}

Using the technique described in \cite{BI03}, the first step is
the construction of an internal model  for  (\ref{fzero}), viewed
as an autonomous system with output (\ref{qzero}). To this end,
consider the sequence of functions recursively defined
 as
 \[\tau_1({\bf z})={\bf q_0}({\bf z})\,,\quad \ldots , \quad \tau_{i+1}({\bf z}) =
 {\partial \tau_i \over \partial {\bf z}}{\bf f}_0({\bf z})
\]
for $i=1, \ldots, d-1$, with $d$ as introduced in assumption (A3),
and consider the map
 \[ \ba{rcccl}
 \tau &:& P\times W\times \Real^n &\to & \Real^d\\[1mm]
&&(\mu,w,z) &\mapsto& {\rm
 col}(\tau_1({\bf z}), \tau_2({\bf z}), \ldots, \tau_d({\bf z}))\,.\ea\] If
 $k$, the degree of continuous differentiability of the functions
 in (\ref{sys1}), is large enough, the map  $\tau$ is well
 defined and $C^1$. In particular $\tau({\mathcal A}_0)$, the image of
 ${\mathcal A}_0$ under $\tau$\, is a {\em compact}
 subset of $\Real^d$, because ${\mathcal A}_0$ is a compact subset of
 $P\times W\times \Real^n$.

 Let $ \varphi_{\rm c}: \Real^d \to \Real$ be any locally Lipschitz
 function of compact support which agrees  on $\tau({\mathcal A}_0)$ with
 the function $ \varphi$ defined in (A3), i.e. a function such that, for
 some compact superset ${\mathcal S}$ of $\tau({\mathcal A}_0)$ satisfies
 \[\ba{rcll}  \varphi_{\rm c}(\eta) &=& 0 & \mbox{for all $\eta \not\in {\mathcal S}$} \\[1mm]
  \varphi_{\rm c}(\eta) &=& \varphi(\eta) & \mbox{for all $\eta \in
 \tau({\mathcal A}_0)$}.\ea
 \]

 With this in mind, consider the system
 \beeq{\label{observer}
\dot \xi =\Phi_{\rm c}(\xi) + G(u_{\rm ss} - \Gamma \xi)}
 in which
 \[
 \Phi_{\rm c}(\xi) = \qmx{\xi_2 \cr \xi_3 \cr \cdots  \cr \xi_d \cr
 - \varphi_{\rm c}(\xi_1,\xi_2,\ldots,\xi_d)\cr},\qquad G = \qmx{ \kappa c_{d-1}\cr
 \kappa^2c_{d-2}
 \cr \cdots  \cr \kappa^{d-1}c_1\cr \kappa^dc_o\cr},
\qquad
\Gamma = \qmx{1 & 0 & \cdots & 0\cr} \,,\]
 the $c_i$'s are such that the polynomial $\lambda^d +
c_0 \lambda^{d-1}+ \cdots + c_{d-1}=0$ is Hurwitz and $\kappa$ is
a  positive number. As shown in \cite{BI03}, if $\kappa$ is large
enough, the state $\xi(t)$ of (\ref{observer}) asymptotically
tracks $\tau({\bf z}(t))$, in which ${\bf z}(t)$ is the state of
system (\ref{fzero}). Therefore $\Gamma \xi(t)$ asymptotically
reproduces its output (\ref{qzero}), i.e. the steady state control
$u_{\rm ss}(t)$. As a matter of fact, the following result holds.

\begin{lemma}\label{LM1} Suppose assumptions (A1) and (A2) hold.
Consider the triangular system
\beeq{\label{observer1}
\ba{rcl} \dot {\bf z} &=& {\bf f}_0({\bf z})\\[1mm]
\dot \xi &=& \Phi_c(\xi) + G({\bf q}_0({\bf z})-\Gamma \xi)\,.\ea}
Let the initial conditions for ${\bf z}$ range in the set ${\bf
Z}$ and  { let $\Xi$ be an  arbitrarily large compact set of
initial condition for $\xi$}. There is a number $\kappa^\ast$ such
that, if $\kappa \ge \kappa^\ast$, the trajectories of
(\ref{observer1}) are bounded and
\[
{\rm graph}(\tau|_{{\mathcal A}_0}) = \omega({\bf Z}\times \Xi)\,.
\]
In particular ${\rm graph}(\tau|_{{\mathcal A}_0})$ is a compact invariant set
which uniformly attracts ${\bf Z}\times \Xi$. Moreover ${\rm
graph}(\tau|_{{\mathcal A}_0})$ is also locally exponentially attractive.
\end{lemma}

{ Note that, as a consequence of assumption (A1), of the fact that
$\tau(\cdot)$ is a continuous vector field and of the fact that
the compact set $\Xi$ can be taken arbitrarily large, it is
possible to assume, without loss of generality, the existence of a
positive $d_2$ such that
\[
 {\bf z} \in P \times W \times \Real^n\,,\quad
 \xi \in \Real^d\,,
 \quad
 |({\bf z},\xi)|_{\omega({\bf Z} \times \Xi)} \leq d_2
 \qquad \Rightarrow \qquad
 ({\bf z},\xi) \in {\bf Z} \times \Xi\,.
\]
}

 In view of Lemma \ref{LM1}, it would be natural -- if the
true error variable $e$ were available for feedback purposes -- to
choose for (\ref{sys1}) a control of the form
 \beeq{\label{controllo0}\ba{rcl} \dot \xi &=& \Phi_{\rm c}(\xi) - Gke \\
 u &=& \Gamma\xi -ke\,, \ea}
 with $k$ a large number. This control, in fact, would solve the
 problem of output regulation (see \cite{BI03}). The true error
 $e$ not being available, we choose instead
 \beeq{ \label{hatedef}
 \hat e = y - y_{\rm r}(w_{\rm d})
 }
 and the controller accordingly as
  \beeq{\label{controllo1}\ba{rcl} \dot \xi &=& \Phi_{\rm c}(\xi) - G k\hat e\\
 u &=& \Gamma\xi -k\hat e\,. \ea}

 The result which will be proven next is that there exists
 $k^\ast>0$ such that if $k\geq k^\ast$ the regulator designed
 above solves the problem in question (provided that $N$ and $T$
 satisfy the condition of Proposition \ref{PR1}).

\section{Main results}

\subsection{Trajectories of the closed loop system are bounded}

To prove that the proposed regulator solves the problem, we show
first of all the trajectories of the controlled system, namely
those of the system \beeq{\label{closedloop} \ba{rcll} \dot w_{\rm
d} &=& s(w_{\rm d})&  \qquad w_{\rm d}(kT) = w_{\rm d}(kT^-) +
w_{\rm q}(k) \dst {L(k)
 \over 2 N}\\[2mm]\dot z &=& f(z,y,\mu)\\[2mm]\dot y &=& q(z,y,\mu) + \Gamma \xi - k(y-
 y_{\rm r}(w_{\rm d}))\\[2mm]
 \dot \xi &=& \Phi_{\rm c}(\xi) - Gk(y-y_{\rm r}(w_{\rm d}))
 \ea} are bounded. To study trajectories of (\ref{closedloop}) it is convenient to
replace the coordinate $y$ by
\[
\hat  e= y - y_{\rm r}(w_{\rm d})\] to obtain the system
\beeq{\label{closedloopdue}\ba{rcl} \dot w_{\rm d} &=&
s(w_{\rm d})\\[1mm]\dot z &=& f(z,\hat e+y_{\rm r}(w_{\rm d}),\mu)\\[1mm]
 \dot \xi &=& \Phi_{\rm c}(\xi)  -Gk \hat e\\[1mm]\dot
 {\hat e}
 &=& q(z,  \hat e+y_{\rm r}(w_{\rm d}),\mu) -L_s y_{\rm r}(w_{\rm d})+ \Gamma \xi -k
 \hat e\,.
 \ea
} This system can be further simplified by changing the state
variable $\xi$ into $\tilde \xi = \xi - G\hat e $ and setting $p =
{\rm col}(\mu,w_{\rm d},z,\tilde \xi)$, so as to obtain a system
of the form
\beeq{\label{short}
\ba{rcl}
\dot p &=& F_0(p) + F_1(p,\hat e)\hat e\\
\dot{\hat e} &=& H_0(p) + H_1(p, \hat e)\hat e - k\hat e\,,\ea}
in which
\[
F_0(p) = \qmx{0 \cr s(w_{\rm d})\cr f(z, y_{\rm r}(w_{\rm d}),\mu)\cr
\Phi_c(\tilde \xi)+G( -q(z,y_{\rm r}(w_{\rm d}),\mu) + L_s y_{\rm r}(w_{\rm
d})-\Gamma \tilde \xi)\cr}\]

\[ H_0(p)=
q(z,y_{\rm r}(w_{\rm d}),\mu) - L_s y_{\rm r}(w_{\rm d}) +\Gamma
\tilde \xi\] and $F_1(p,\hat e)$, $H_1(p,\hat e)$ are suitable
continuous functions.

{
 With this notation in mind, we state the next proposition
 which claims that a large value of $k$ succeeds in rendering bounded the
 trajectories of the switched nonlinear system (\ref{short})
 provided that the sampling interval $T$ is sufficiently large.}

\begin{proposition}\label{PR2}
  Consider system (\ref{closedloop}) with initial conditions in
  $P\times W \times Z \times Y \times \Xi$. Suppose assumptions (A0)-(A3) hold.
  Let $\kappa$ be chosen as indicated in Lemma \ref{LM1}.
 { Then there exist $T^\ast>0$ and $k^\ast>0$ such that for all sampling
  intervals $T > T^\ast$ and all  $ k \geq k^\ast$ the trajectories
  are bounded in positive time.}
\end{proposition}

\begin{proof} System (\ref{short}) is apparently identical to the closed-loop system already
studied in \cite{BI03}, but with the exception  that, every $T$
units of time, $w_{\rm d}$ is being ``reset" as
\[
w_{\rm d}(kT) = w_{\rm d}(kT^-) + w_{\rm q}(k) \dst {L(k)
 \over 2 N}
\] and
$\hat e$ is being ``reset" as
\[\hat e(kT) = \hat e(kT^-) - y_{\rm r}(w_{\rm d}(kT^-) + w_{\rm q}(k) \dst {L(k)
 \over 2 N})+ y_{\rm r}(w_{\rm d}(kT^-))\,.
\]

With this in mind, we proceed to study the behavior of
(\ref{short}) on the time interval $[0,T)$. Observe that system
\beeq{\label{psys}
\dot p = F_0(p)}
coincides with system  (\ref{observer1}), the only difference
being that the component $w$ of ${\bf z}$ is now written as
$w_{\rm d}$ and $\xi$ is now written as $\tilde \xi$. Thus, from
Lemma \ref{LM1}, it can be asserted that in this system all
trajectories with initial conditions in ${\bf Z}\times
\Xi$ are attracted by the compact invariant set
\beeq{\label{limitset}
{\mathcal A}=\{({\bf z},\xi) \in {\mathcal A}_0\times \Real^d:
\xi=\tau({\bf z})\}\,.}
Moreover, by construction, the function $H_0(p)$ vanishes on this
set.

Let ${\mathcal D}$ denote the domain of attraction of ${\mathcal
A}$. Then, as shown for instance in \cite{BIP}, for any given
arbitrarily small $\epsilon$, it is possible to claim the
existence of a continuous function $V: {\mathcal D} \to \Real$
having the following properties:

\medskip\noindent
(a) $V(p)=0$ if $p\in {\mathcal A}$ and $V(p)>0$ everywhere else,

\medskip\noindent
(b) $V(\cdot)$ is proper on ${\mathcal D}$,
\medskip\noindent

 {
\medskip\noindent
(c) for some positive $b$
\[
 V(p) \leq b \qquad \Rightarrow \qquad |p|_{\cal A} \leq
 \epsilon\,,
\]

\medskip\noindent
(d) for some positive $g < b$, $V(\cdot)$ is locally Lipschitz on
the set ${\mathcal D}_g=\{p\in {\mathcal D}: V(p)> g\}$,

\medskip\noindent
(e) for any $p \in {\mathcal D}_g$\,,
\[ \limsup_{h\to 0^+} {1\over
h}[V(\phi(h,p))-V(p)] \le -1 \] in which $\phi(t,p)$ denotes the
flow of (\ref{psys}).
}

\medskip\noindent
{
 Let now $\epsilon < d_2$ (with $d_2$ defined after
Lemma \ref{LM1}), pick the function $V$ accordingly and note that
by construction \beeq{ \label{pb}
 V(p) \leq b \qquad \Rightarrow \qquad |p|_{\cal A} \leq
 \epsilon \leq  d_2 \qquad \Rightarrow \qquad
 p \in {\bf Z}\times \Xi\,.
}
 Moreover pick numbers $a$ and $b_1$, with $a>b$  such that ${\bf Z}\times
\Xi\subset V^{-1}([0,a])$ (which is possible as $V$ is proper on
$\cal D$), and $g < b_1 < b$, and set
 \[
 {\mathcal S} =\{ p \in {\mathcal D} \quad : \quad b_1 \leq V(p) \leq a+1
 \}\,.
 \]
}

{
 Let $c$ and $\Delta$ be positive numbers such that
 \[
 |y-y_{\rm r}(w)|\leq c \qquad \mbox{for all } y\in Y \mbox{ and }
 w \in W_0\]
 and
 \[
 |y_{\rm r}(w)-y_{\rm r}(w_{\rm d})| \leq \Delta \qquad \mbox{for all } w_{\rm d}\in W \mbox{ and }
 w \in W_0
 \]
 and note that
 \[
 \ba{rcl}
 |\hat e(0)| = |y(0)-y_{\rm r}(w_{\rm d}(0))| &=&
 |y(0)-y_{\rm r}(w(0)) + y_{\rm r}(w(0)) - y_{\rm r}(w_{\rm d}(0))|\\[1mm]
 & \leq& |y(0)-y_{\rm r}(w(0))| + |y_{\rm r}(w(0))-y_{\rm r}(w_{\rm d}(0))| \leq c + \Delta\,.
\ea
 \]

 Finally, let $\bar f$ and  $\bar h$ be defined as
 \[
 \bar f := \max_{{p \in V^{-1}([0,a+1])\atop |\hat e|\leq {c+\Delta+1}}}
 |F_1(p,\hat e)|\qquad
 \bar h := \max_{{p \in V^{-1}([0,a+1])\atop |\hat e|\leq {c+\Delta+1}}}
 |H_0(p)+H_1(p,\hat e)\hat e|\,.
 \]
}

 {\em Claim 1: There exists $T^\ast>0$ and $k^\ast>0$ such that
 for any $T\geq T^\ast$, any $k\ge k^\ast$ and any $\ell\geq 0$}
 \beeq{\label{claim1}
 \ba{c}|\hat
 e(\ell T )|\le c+\Delta \\[1mm]
 p(\ell T)\in {\bf Z} \times \Xi
 \ea\qquad
 \Rightarrow \qquad \ba{c}|\hat
 e((\ell +1)T^-)|
 \le c\\[1mm]
 p((\ell + 1)T^-)\in
 {\bf Z} \times \Xi \ea\;\;.\qquad \triangleleft}

 To prove the claim, we
 derive two basic inequalities. The first one is about the
 dynamics of $\hat e$. Looking at the bottom equation of
 (\ref{short}) and bearing in mind the definition of $\bar h$,
 standard arguments can be used to claim that -- if on some time
 interval $[t_0, (\ell + 1)T)$, $t_0 \geq \ell T$,
 the state of (\ref{short}) satisfies $|\hat
 e(t)|\le c+\Delta+1$ and $p(t)\in V^{-1}([0,a+1])$ -- then
 \beeq{ \label{bx}
  |\hat e(t)|
   \leq e^{- k (t-t_0)}|\hat e(t_0)| + {\bar h \over k}\qquad\qquad
   \forall t\in [t_0, (\ell + 1)T)\,.
 }

 The second inequality is about the dynamics of $p$. Looking at the top equation of
 (\ref{short}) and taking the Dini's derivative of $V(p(t))$ we see
 that -- if on some time
 interval $[t_0, (\ell + 1)T)$, with $t_0 \geq \ell T$, the state of (\ref{short})
 satisfies $|\hat e(t)|\le c+\Delta+1$ and $p(t)\in {\mathcal S}$ -- then
 \beeq{\label{Dpb}
 D^+V(p(t)) \leq - 1 + L_V \bar f |\hat e(t)| \leq - 1 + L_V \bar
 f
 (e^{- k (t-t_0)}|\hat e(t_0)| + {\bar h \over k})\qquad\qquad
   \forall t\in [t_0, (\ell + 1)T)\,,
 }
 in which $L_V$ is such that $|V(p) - V(q)| \leq L_V |p-q|$
 $\forall$ $p,q \in \mathcal S$. From this, using the appropriate comparison
 lemma, we obtain
 \beeq{ \label{bp}
 V(p(t)) \leq V(p(t_0)) - (1 - {L_V \bar f\, \bar h \over k})(t-t_0) +
 {L_V\bar f \,|\hat e(t_0)| \over k}\qquad\qquad
   \forall t\in [t_0, (\ell + 1)T)\,.
 }

In order to be able to prove Claim 1, we need this auxiliary
result.

\medskip
{\em Claim 2: there exists $k^{\ast\ast}$ such that for any $k\ge
k^{\ast\ast}$}, $T \geq T^\ast$ and $\ell \geq 0$
 \beeq{\label{claim11}
 \ba{c}|\hat
  e(\ell T)|\le c+\Delta\\[1mm]
  p(\ell T)\in V^{-1}([0,a])\ea \qquad \Rightarrow \qquad
  \ba{c}|\hat e(t)|\leq c+\Delta+1 \\[1mm]
   p(t) \in V^{-1}([0,a+1])\ea \qquad \forall \, t \in [\ell T, (\ell + 1)T)\,.\qquad\triangleleft}
 This claim can be proved by contradiction. As a matter of fact
 suppose that (\ref{claim11}) is not true, namely that there exists a time
 $T' \in [\ell T, (\ell + 1)T) $
 such that either $|\hat e(T')| > c+\Delta+1$ or $V(T') > a+1$ (with a
 mild abuse of notation we write $V(\cdot)$ for $V(p(\cdot))$).
 By continuity of the trajectories with respect to time, this means that
 there exists a time $T'' \geq \ell T$, $T'' \leq T'$, such that either
  \beeq{\label{c1}
 |\hat e(t)|< c+\Delta+1 \quad t \in [\ell T,T'')\, \quad
 |\hat e(T'')| = c+\Delta+1 \qquad \mbox{and} \qquad V(t)
 \leq a+1 \quad t \in [\ell T,T'']
 }
 or
 \beeq{ \label{c2}
 V(t)< a+1 \quad t \in [\ell T ,T'')\, \quad V(T'') = a+1 \qquad \mbox{and} \qquad
 |\hat e(t)|
 \leq c+\Delta+1 \quad t \in [\ell T,T'']\,.
 }
But if (\ref{c1}) were true, by (\ref{bx}) taking $t_0 = \ell T$
and $k \geq 2 \bar h$,
 we would have that
 \[
 |\hat e(T'')|\leq e^{-k (T''-\ell T)}|\hat e(\ell T)| +
 {\bar h \over k} \leq c+\Delta + {\bar h \over
 k} \leq c +\Delta+ {1\over 2}
 \]
 which contradicts $|\hat e(T'')|=c+\Delta+1$. A similar contradiction would be obtained
  if (\ref{c2}) were true. In fact, let $t_0$, $\ell T \leq
  t_0<T''$,
  be any time such
  that $V(t_0)=a$ and $V(t)\ge a$ for all $t\in [t_0,T'']$. Using
  (\ref{bp}) and taking $k \geq \max \{L_V \bar f \,\bar h, \;2 L_V\bar f (c+1)\}$,
 we have that
 \[
 V(p(T'')) \leq a - (1 - {L_V \bar f\, \bar h \over k})(T''-t_0) + {L_V\bar f (c+1) \over
 k} \leq a + {1 \over 2}
 \]
 which contradicts $V(T'') = a+1$. From this, Claim 2 follows by taking
 \[
 k^{\ast\ast} = \max \{2 \bar h, \; L_V \bar f\, \bar h, \; 2 L_V\bar f (c+\Delta+1) \}\,.
 \]
 $\diamond${\em (End proof Claim 2)}\\[2mm]

 \noindent Having proven that  $|\hat e(t)| \leq c+\Delta+1$ and $p(t) \in V^{-1}([0,a+1])$
 on the entire time interval $[\ell T, (\ell +1)T)$,
 the desired Claim 1  follows again from arguments based on the
 bounds (\ref{bx}) and (\ref{bp}). To show that a large value of $k$ succeeds in
recovering $|\hat e((\ell+1)T^-)|\leq c$, just set
\[
 k \geq \max\{{2 \bar h \over  c}, \; {1 \over T}\ln {2(c+\Delta)  \over c} \}
\]
in the estimate (\ref{bx}) evaluated with $t_0 = \ell T$.

To show that a large value of $k$ succeeds in recovering
$p((\ell+1)T^-) \in {\bf Z} \times \Xi$ we prove that if
 \[
 T \geq T^\ast = 2(a-b) + 1
 \]
  then $V(p((\ell + 1)T^-)) \leq b $ which, by (\ref{pb}), implies
  $p((\ell + 1)T^-) \in {\bf Z} \times \Xi$. To this end, we
  distinguish two cases. If $p(t)\in {\mathcal S}$ for all $t\in [\ell T,
  (\ell+1)T)$, just set
 \[
 k \geq \max\{ 2 L_V \bar f (c+\Delta+1)\,, \;  2 L_V \bar f\, \bar h\}
 \]
in the estimate (\ref{bp}) with $t_0 = \ell T$ to obtain $V((\ell +1)T^-)\le
b$. The condition $p(t)\in {\mathcal S}$ for all $t\in [\ell T, (\ell +1)T)$
can be violated if and only if there are times in $[\ell T, (\ell +1)T)$ at
which $V(t) < b_1$ (recall the definition of ${\mathcal S}$ and the fact that
$V(t)\le a+1$ for all $t\in [\ell T, (\ell +1)T)$). Let this be the case. If
$V((\ell +1)T^-)< b_1 < b$, the claim is trivially true. If not, let $T' \in
[\ell T, (\ell +1)T)$ be such that $V(T')=b_1$ and $b_1\le V(t)\le a+1$ for all
$t\in [T', (\ell +1)T)$. On this time interval one can still use (\ref{bp})
with $t_0 = T'$ to conclude that if $k \geq {L_V \bar f (c+\Delta+1)/(b-b_1)}$
then $V(t - T') \leq b$ for all $t \in [T', (\ell +1)T)$ from which the Claim 1
follows. $\diamond${\em (End proof of Claim 1)}\\[2mm]

\noindent With this result at hand, Proposition 2 can be easily
proved by subsequent applications of Claim 2.
 To this end note that $p(0) \in
{\bf Z} \times \Xi$ and, by definition of $c$ and $\Delta$, $|\hat e(0)|\leq c
\leq c+\Delta$, from which Claim 1 evaluated for $\ell=0$ yields $|\hat e(T^-)|
\leq c$ and $p(T^-) \in {\bf Z} \times \Xi$. Suppose now that, for some $\ell >
0$, $|\hat e(\ell T^-)| \leq c$ and $p(\ell T^-) \in {\bf Z} \times \Xi$.
 Then, after the switch,
\[
\ba{rcl}
 |\hat e(\ell T)| &=& |y(\ell T) - y_{\rm r}(w_{\rm d}(\ell T))| =
 |y(\ell T^-) - y_{\rm r}(w_{\rm d}(\ell T))|\\[1mm]
 &=& |y(\ell T^-) - y_{\rm r}(w_{\rm d}(\ell T^-))
 +y_{\rm r}(w_{\rm d}(\ell T^-)) - y_{\rm r}(w_{\rm d}(\ell T)|\\[1mm]
 &\leq& |\hat e(\ell T^-)| + |y_{\rm r}(w_{\rm d}(\ell T^-)) - y_{\rm r}(w_{\rm d}(\ell T))|
 \leq |\hat e(\ell T^-)| + \Delta \leq c+\Delta
 \ea
\]
 and, by item (i) of Proposition \ref{PR1}, $w_{\rm d}(\ell T) \in W$. The
 latter, by bearing in mind the definition of the set $\bf Z$ and
 the definition of $p$, yields $p(\ell T) \in {\bf Z} \times \Xi$
 from which the result of Proposition \ref{PR2} follows.
 $\triangleleft$
\end{proof}

Proposition 2 shows that trajectories of the controlled system remain bounded
if the time interval $T$ exceeds a minimum number $T^\ast$ (minimal
``dwell-time'') which depends on the parameters of the controlled system and on
the sets of initial conditions. This, however, may be in contrast with relation
(\ref{Nbconst}) which, bearing in mind that $M(\cdot)$ is an increasing
function, requires that the sampling interval $T$ is small enough. In the case
the dwell-time $T^\ast$ is not compatible with (\ref{Nbconst}), a simple
modification of the decoder structure helps solving the problem.
 Let $N_b$ be given, let $\bar T$ denote the minimal value of
 $T$ compatible with (\ref{Nbconst}) and let $\ell$ be any
 positive integer such that
 \[
 \ell \bar T \geq T^\ast\,.
 \]
 Consider now a ``second level'' decoder dynamics defined as
\beeq{\label{sldecdyn}
 \dot w_{\rm d}' = s(w_{\rm d}') \qquad w_{\rm d}'(0) = w_{\rm
 d}(0)
 }
 whose state $w_{\rm d}'$ is periodically reset, every $\ell \bar T$ units of time,
 to the value of the ``first level'' decoder dynamics
 (\ref{decdyn}), that is as
 \[
 w_{\rm d}'(k \ell \bar T) = w_{\rm d}(k \ell \bar T) \qquad \mbox{for all }k\geq
 0\,.
 \]
  In other words $w_{\rm d}'$ provides an under-sampled version of the
 first-level decoder dynamics (\ref{decdyn}) with the under-sampling period
 $\ell \bar T$ such that the constraint on the minimal dwell-time is respected.

 The same arguments used to prove Proposition 1 show that $w_{\rm d}'(t) \in W$
 for all $t \geq 0$ and furthermore
 \[
 \lim_{t \rightarrow \infty} |w_{\rm d}'(t) - w(t)| =0
 \]
 with uniform convergence rate.

 Consider now the regulator (\ref{controllo1}) in which $\hat e$ is defined
 as in (\ref{hatedef})
 but with $w_{\rm d}$ replaced by $w_{\rm d}'$. It is easy to realize
 that all the analysis carried out in this section can now be repeated
 by replacing the decoder
 dynamics (\ref{decdyn}) by (\ref{sldecdyn}), the decoder
 variable $w_{\rm d}$ by $w_{\rm d}'$, and the time interval $T$
 by $\ell \bar T$. In particular by Proposition 2 and  by the fact that
 $\ell \bar T \geq T^\ast$ we conclude that the trajectories of
 the controlled system are bounded if $k$ is chosen sufficiently
 large.
  This and the fact that $w_{\rm d}'(t)$ converges asymptotically
  to $w(t)$ with uniform convergence rate are the crucial properties
  needed to conclude that the proposed regulator solves the
  problem  in question as precisely stated and proved in the next subsection.

\subsection{The tracking error converges to zero}

 To prove that the tracking error converges to zero, it is useful
 to observe that, if the coordinate $y$ of (\ref{closedloop}) is
 replaced by
 \[
 e = y - y_{\rm r}(w)\]
 the system in question can be also rewritten as
 \beeq{\label{esys}\ba{rcl}
\dot w &=&
s(w)\\[1mm]\dot z &=& f(z,e+y_{\rm r}(w),\mu)\\[1mm]
 \dot \xi &=& \Phi_{\rm c}(\xi) + G(-k e)+G(-k\tilde e)\\[1mm]\dot
 e
 &=& q(z,  e+y_{\rm r}(w),\mu) -L_s y_{\rm r}(w)+ \Gamma \xi -k
 e-k\tilde e
 \ea
} having set \[ \tilde e = \hat e - e\,.
\]

The same change of variables used to put (\ref{closedloopdue}) in
the form (\ref{short}) yields now a system of the form
\beeq{\label{shortpert}
\ba{rcl}
\dot p &=& F_0(p) + F_1(p,e)e\\
\dot e &=& H_0(p) + H_1(p, e)e - ke-k\tilde e\,,\ea}
in which $p={\rm col}(\mu,w,z,\tilde \xi)$ and $F_0(p), F_1(p,e),
H_0(p), H_1(p, e)$ are the same as in (\ref{short}). This system
can be viewed as a ``perturbed" version of system
\beeq{\label{short1}
\ba{rcl}
\dot p &=& F_0(p) + F_1(p,e)e\\
\dot e &=& H_0(p) + H_1(p, e)e - ke\ea}
whose asymptotic properties have been investigated in \cite{BIP}.

The following result is a minor enhancement of the main result of
\cite{BIP}. Let $V(\cdot)$ be the positive definite function
introduced in the proof of Lemma \ref{LM1}, set ${\mathcal P}=\{p:
V(p)\le a\}$ with $a$ chosen so that $P\times W\times Z \times \Xi
\subset {\mathcal P}$, and set $E=\{e: |e|\le c\}$. Moreover, let
${\mathcal A}$ be the set defined by (\ref{limitset}).

\begin{lemma}\label{LM2} Consider system
(\ref{short1}) in which $F_0(p), F_1(p,e), H_0(p), H_1(p, e)$ are defined as
before and initial conditions are taken in ${\mathcal P}\times E$. Suppose
assumptions (A0)-(A3) hold. Let $\kappa$ be chosen as indicated in Lemma
\ref{LM1}. Then there is $k^\ast$ such that, if $k>k^\ast$, the following
holds:

\medskip\noindent {\rm (i)} the positive orbit of ${\mathcal
P\times E}$ under the flow of (\ref{short1}) is bounded and
\[
\lim_{t\to \infty} |p(t)|_{\mathcal A} = 0, \qquad \lim_{t\to
\infty}e(t) =0\,.
\]

\medskip\noindent {\rm (ii)} for any $\varepsilon >0$, there exist
numbers $\delta_1>0$ and $\delta_2>0$ such that, if
$|p_0|_{\mathcal A}\le \delta_1$ and $|e_0|\le \delta_1$, for any
continuous function $u(t)$ satisfying $|u(t)|\le \delta_2$ for all
$t\ge 0$, the solution $p(t),e(t)$ of the perturbed system
\beeq{\label{short1pert} \ba{rcl}
\dot p &=& F_0(p) + F_1(p,e)e\\
\dot e &=& H_0(p) + H_1(p, e)e - ke+u(t)\ea}
with initial conditions $p(0)=p_0$ and $e(0)=e_0$ satisfies
\[
 |p(t)|_{\mathcal A} \le \varepsilon,\qquad |e(t)|\le \varepsilon,
\qquad \forall t\ge 0\,.
\]
\end{lemma}

\begin{proof}
Item (i) has already been proven in \cite{BIP}. The proof of item (ii) consists
in a minor modification of the arguments used in \cite{BIP} to prove item (i).
Consider the Locally Lipschitz Lyapunov function $U(p)$ defined in \cite{BIP},
which -- for some $a_1<1$, $\lambda>0$ and $\bar L>0$ -- satisfies
\[
a_1|p|_{\mathcal A} \le U(p) \le |p|_{\mathcal A}\,,\]  and
\[
D^{\,+}U(p(t)) \le -\lambda U(p(t)) + \bar L \bar f |e(t)|
\]
along the integral curve $p(t)$, so long as  $p(t)$ remains
sufficiently close to ${\mathcal A}$ (see (24) in \cite{BIP}).

Consider now, for (\ref{short1pert}), the candidate Lyapunov
function
\[
W(p,e) = {1\over 2}(U^2(p) + e^2)
\]
which trivially satisfies
\[
{a_1 \over 2}(|p|^2_{\mathcal A} + e^2) \le W(p,e) \le {1 \over
2}(|p|^2_{\mathcal A} + e^2)\,.
\]

Taking its Dini derivative along the trajectories of
(\ref{short1pert}), we can obtain an estimate of the form (we omit
the argument $t$ for convenience)
\[
D^{\,+}W(p,e) \le -\lambda U^2(p) + \bar L \bar f U(p) |e| + \beta
|p|_{\mathcal A}|e| - (k-\bar k) |e|^2 + |u||e|\]
 so long as $p(t)$ remains sufficiently close to ${\mathcal A}$, in which $\beta$
 and $\bar k$ are fixed positive numbers. Using the
above estimates for $U(p)$, it is easy to deduce that
\[
D^{\,+}W(p,e) \le \qmx{U(p) & |e|\cr}\qmx{-\lambda & {1\over
2}(\bar L \bar f + {\beta\over a_1}) \cr {1\over 2}(\bar L \bar f
+ {\beta\over a_1}) & -(k-\bar k- {1\over 2})}\qmx{U(p) \cr
|e|\cr} + {1\over 2}|u|^2\,.
\]

Clearly, there is a value $k^\ast>0$ and a number $a>0$ such that,
if $k\ge k^\ast$
\[
\qmx{U(p) & |e|\cr}\qmx{-\lambda & {1\over 2}(\bar L \bar f +
{\beta\over a_1}) \cr {1\over 2}(\bar L \bar f + {\beta\over a_1})
& -(k-\bar k- {1\over 2})}\qmx{U(p) \cr |e|\cr} \le - {a\over
2}(U^2(p)+e^2) \le - aW(p,e)\,,
\]
and this yields, using the appropriate comparison lemma,
\[
W(p(t),e(t)) \le - e^{-at}W(p_0,e_0) + {1\over 2a}\max_{\tau\in
[0,t]}|u^2(t)|\,,
\]
for all $t\ge 0$. From this, the result follows by standard
arguments. $\triangleleft$
\end{proof}

We are now ready to prove the main result of the paper.

\begin{proposition} Consider system (\ref{closedloop}) with initial conditions in $P\times W \times Z \times \Xi
 \times Y$.
Suppose assumptions (A0)-(A3) hold. Let $\kappa$ be chosen as indicated in
Lemma \ref{LM1} and $k$ as indicated in Proposition \ref{PR2} and in Lemma
\ref{LM2}. Then
\[
\lim_{t \to \infty}e(t) = 0\,.
\]
\end{proposition}

\begin{proof}
As already mentioned, system (\ref{closedloop}) can be written in
the form (\ref{shortpert}). Moreover, if the initial condition of
(\ref{closedloop}) is taken in the set $P\times W \times Z \times
\Xi \times Y$, the corresponding initial condition of (\ref{shortpert}) is in ${\mathcal P}\times E$.
 In view of the result of
Lemma \ref{LM2}, item (ii), the result is proven if we are able to
show that, given any pair of numbers $\delta_1$ and $\delta_2$,
there is a time $\bar t$ such that
 \beeq{\label{delta1} |p(\bar
t)|_{\mathcal A} \le \delta_1,\qquad |e(\bar t)| \le \delta_1 }
and \beeq{\label{delta2} |k\tilde e(t)| \le \delta_2, \qquad
\forall t\ge \bar t\,. }

To this end, recall that by definition
\[
\tilde e = \hat e - e = -y_{\rm r}(w_{\rm d}) + y_{\rm r}(w)\,.
\]
Hence, since $y_{\rm r}(\cdot)$ is continuous, in view of Proposition \ref{PR1}
we have that $\lim_{t\to \infty}\tilde e(t) =0$. As a consequence, there is a
time $t^\ast$, dependent on $\delta_2$, such that (\ref{delta2}) is fulfilled
for all $\bar t\ge t^\ast$ (note that the coefficient $k$ - which is possibly a
large number - is now fixed). Thus the only critical issue is to make sure that
(\ref{delta1}) holds for some $\bar t\ge t^\ast$. To this end, one can use the
following argument, suggested in \cite{Sontag}.

Consider a system \beeq{\label{s1} \dot x = f(x)+u(t)} in which
$f(\cdot)$ is locally Lipschitz and $u(t)$ is a
piecewise-continuous function. Let $x(t,t_0,x_0,u)$ denote the
integral curve passing through $x_0$ at time $t=t_0$. Suppose
$u(t)$ satisfies
\[
\lim_{t\to \infty}u(t)=0\,,\]
and that, for a  given $x_0$ and a given $T>0$, there is a compact
set $X$ such that
\begin{eqnarray}
x(t,0,x_0,u)\in X & \qquad \forall t\ge 0\,,\label{prop1}\\[1mm]
x(t,\ell T,x(\ell T,0,x_0,u),0)\in X&\qquad \forall t\ge \ell
T,\quad
\forall\, \ell \in \Integer\,.\label{prop2}\end{eqnarray}

\medskip

{\em Claim 3: for any $t_2>0$ and any $\delta>0$, there is a
$\ell^\ast$ such that, for all $\ell\ge \ell^\ast$},
\beeq{\label{closex}
|x(t_2+\ell T,0,x_0,u)- x(t_2+\ell T,\ell T,x(\ell
T,0,x_0,u),0)|\le
\delta\,.\qquad \triangleleft
}

To prove this claim, set $x(t)=x(t,0,x_0,u)$ and $\hat
x(t)=x(t,\ell T,x(\ell T,0,x_0,u),0)$.  For all $t\ge \ell T$, we
have
\[
x(t) = x(\ell T) + \int_{\ell T}^t f(x(s))ds + \int_{\ell T}^t
u(s)ds\] and
\[
\hat x(t) = x(\ell T) + \int_{\ell T}^t f(\hat x(s))ds
\]
Since $f(\cdot)$ is locally Lipschitz and $X$ is compact, there is
$L$ such that $|f(x)-f(y)|\le L|x-y|$ for all $x$ and $y$ in $X$.
Then
\[
|x(t)-\hat x(t)|\le L\int_{\ell T}^t|x(s)-\hat x(s)|ds +
\int_{\ell T}^t |u(s)|ds\,.
\]
Since $u(t)$ converges to 0 as $t\to\infty$, given any
$\varepsilon$ there exists $T_\varepsilon$ such that $|u(t)|
\le
\varepsilon$ for all $t\ge T_\varepsilon$. For any $\ell T\ge T_\varepsilon$ we
have
\[
|x(t)-\hat x(t)|\le L\int_{\ell T}^t|x(s)-\hat x(s)|ds +
\varepsilon(t-\ell T) \,,
\]
from which Gronwall-Bellman's lemma yields
\[
|x(t)-\hat x(t)| \le {\varepsilon \over L}(e^{L(t-\ell T)}-1)
\]
with $t\ge \ell T$, hence \[ |x(\ell T+t_2)-\hat x(\ell T+t_2)|\le
{\varepsilon
\over L}(e^{Lt_2}-1)\,.
\]
Given any $t_2>0$ and $\delta >0$, set $\varepsilon = \delta
L/(e^{Lt_2}-1)$ and choose $\ell^\ast$ so that $\ell^\ast T\ge
T_\varepsilon$. This proves Claim 3.
\medskip

With this result at hand, it is easy to complete the proof of the
Proposition. To this end, we set $x={\rm col}(p,e)$, and identify
system (\ref{short1}) with a system of the form
\beeq{\label{s2}\dot x = f(x)} and system (\ref{shortpert}) with a system of the form
 (\ref{s1}). The assumptions
under which the previous claim holds are satisfied, with $X$ taken
as the set of all $(p,e)$ such that $p\in V^{-1}([0,a+1])$ and
$|e|\le c+1$. In fact, the proof of Proposition \ref{PR2} shows
that trajectories of (\ref{shortpert}) are contained in this set
$X$, i.e. that condition (\ref{prop1}) holds. The same proof also
shows that at each time $\ell T$, where $T$ is the sampling
interval (of the over-sampling period according to the discussion
at the end of the previous subsection), trajectories of
(\ref{shortpert}) are such that $p\in V^{-1}([0,a])$ and $|e|\le
c$. This and the fact that system (\ref{short1}) coincides with
(\ref{shortpert}) with $u =0$  in turn guarantees, by Proposition
\ref{PR2}, that also condition (\ref{prop2}) holds. Since
$x(t,0,x(\ell T,0,x_0,u),0)$ is a solution of (\ref{short1}) with
initial condition in ${\mathcal P}\times E$, we know from Lemma
\ref{LM2}, item (i), that -- given any $\delta_1>0$ - there exists
$t_2$ such that the $p$ and $e$ components of $x(t_2,0,x(\ell
T,0,x_0,u),0)$ satisfy
\[
|p(t_2)|_{\mathcal A} \le {\delta_1 \over 2},\qquad |e(t_2)|\le
{\delta_1 \over 2}\,.
\]

Using now the Claim 3 with $\delta = \delta_1/2$ and the fact that
\[
x(t_2+\ell T,\ell T,x(\ell T,0,x_0,u),0)= x(t_2,0,x(\ell
T,0,x_0,u),0)\,,\] we deduce that the $p$ and $e$ components of
$x(t_2+\ell T,0,x_0,u)$ satisfy, for all $\ell \ge \ell^\ast$
\[
|p(t_2+\ell T)|_{\mathcal A} \le \delta_1,\qquad |e(t_2+\ell
T)|\le \delta_1 \,.
\]
This is what was needed to complete the proof of the Proposition.
$\triangleleft$
\end{proof}

\section{Simulation Results}
 We consider the problem of synchronizing two oscillators located at remote
 places through a constrained communication channel.
 The master oscillator (playing the role of exosystem) is a Van der Pol oscillator
 described by
 \beeq{ \label{exoexa}
 \ba{rcl}
 \dot w_1 &=& w_2 + \epsilon (w_1 + a w_1^3)\\
 \dot w_2 &=& -w_1
\ea
 }
 whose output $y_{\rm ref} = w_2$  must be replied by the output $y$ of a
 remote system of the form
 \beeq{ \label{sysexa}
 \dot y= u\,.
 }
 Simple computations show that, in this specific case, the steady
 state control input $u_{\rm ss}$  coincides with $u_{\rm ss}=
 -w_1$ and the assumption (A3) is satisfied by
 \[
 \ba{rcl}
 \dot \xi_1 &=& \xi_2\\
 \dot \xi_2 &=& f(\xi_1, \xi_2)\\[2mm]
 u_{\rm ss} &=& \xi_1
 \ea
 \]
 where $f(\xi_1,\xi_2) = - \xi_1 - \epsilon(\xi_2 - 3 a \xi_1^2
 \xi_2)$ through the map
 \[
 \tau(w) =\left ( \ba{cc} -w_1 &  -w_2 + \epsilon\, (w_1 - a w_1^3)
 \ea \right )\tr
 \]

\begin{figure}[htb]
\centering
\includegraphics[width=1\textwidth,height=0.5\textwidth]{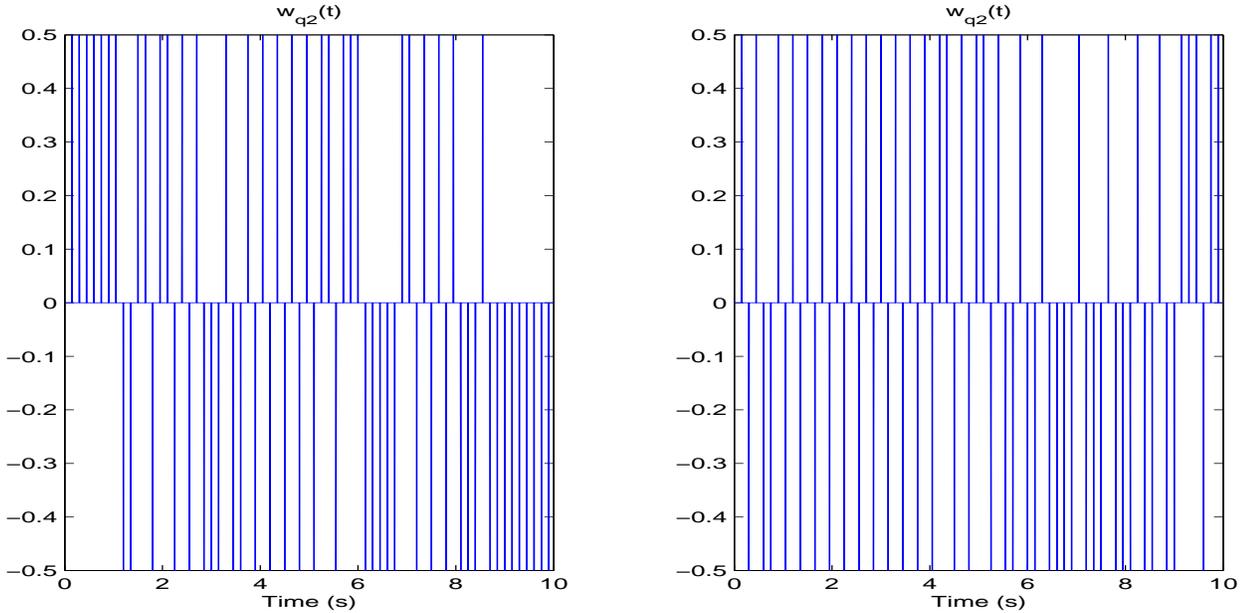}
\caption{\label{1cn1} First control scenario ($N=2$, $T=0.15$ s): behavior of
the encoded variables $w_{q1}(t)$ (left) and $w_{q2}(t)$ (right).}
\end{figure}

\begin{figure}[htb]
\centering
\includegraphics[width=1\textwidth, height=0.5\textwidth]{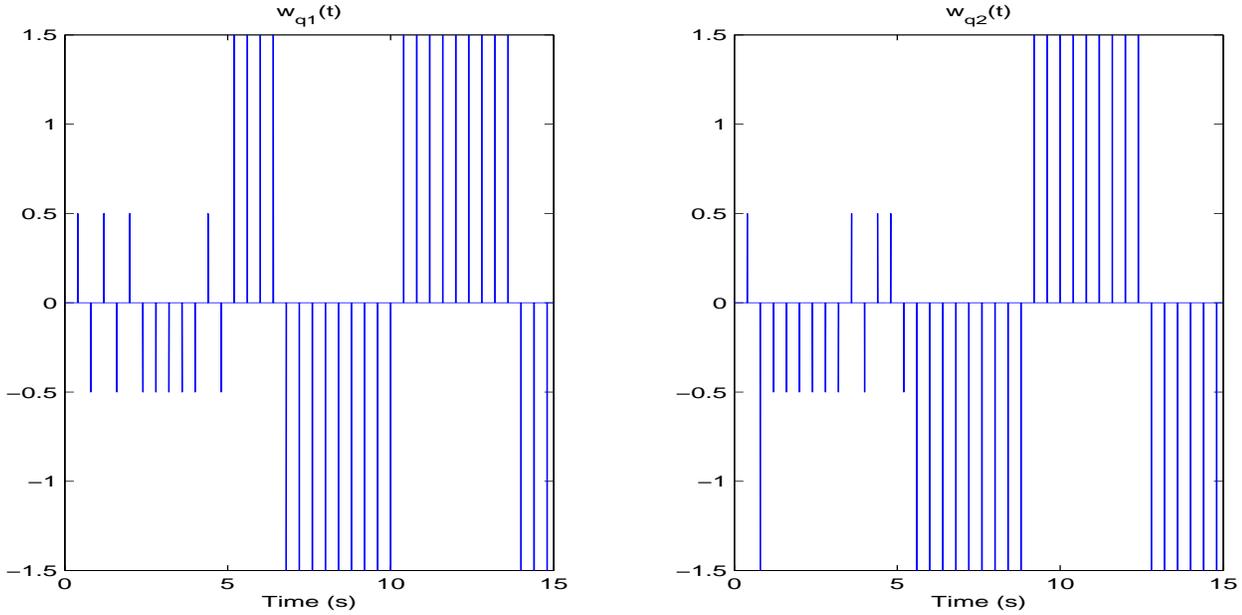}
\caption{\label{2cn1} Second control scenario ($N=4$, $T=0.5$ s): behavior of
the encoded variables $w_{q1}(t)$ (left) and $w_{q2}(t)$ (right).}
\end{figure}

\begin{figure}[htb]
\centering
\includegraphics[width=1\textwidth,height=0.5\textwidth]{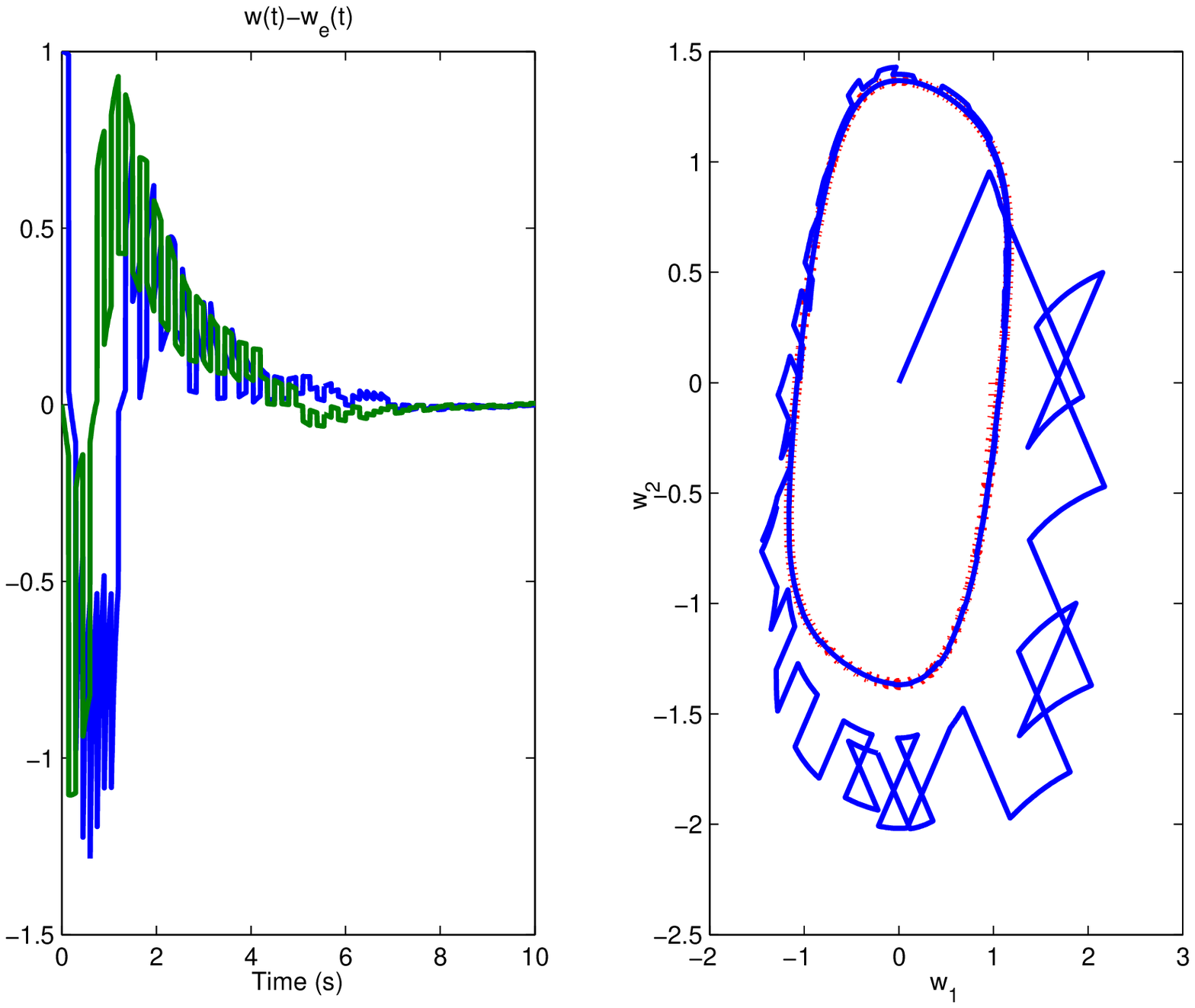}
\caption{\label{1cn2} First control scenario ($N=2$, $T=0.15$ s). Left: time
behavior of $w(t)-w_{\rm e}(t)$ ($w(t)-w_{\rm d}(t)$). Right: phase portrait of
the exosystem (dotted line) and trajectory $(w_{{\rm e}1},w_{{\rm e}2})$ (solid
line).}
\end{figure}

\begin{figure}[htb]
\centering
\includegraphics[width=1\textwidth,height=0.5\textwidth]{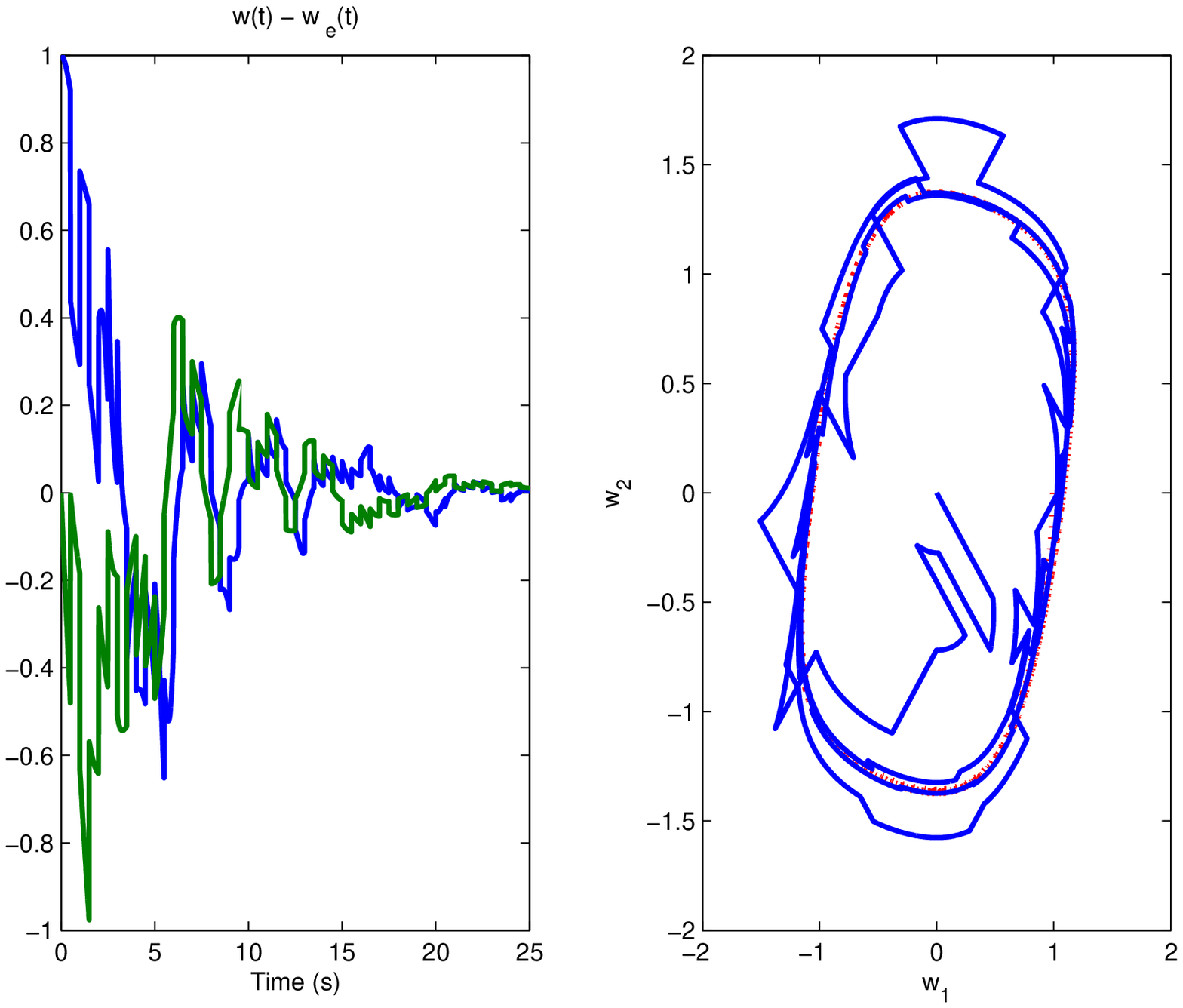}
\caption{\label{2cn2} Second control scenario ($N=4$, $T=0.5$ s). Left: time
behavior of $w(t)-w_{\rm e}(t)$ ($w(t)-w_{\rm d}(t)$). Right: phase portrait of
the exosystem (dotted line) and trajectory $(w_{{\rm e}1},w_{{\rm e}2})$ (solid
line).}
\end{figure}

\begin{figure}[htb]
\centering
\includegraphics[width=1\textwidth,height=0.5\textwidth]{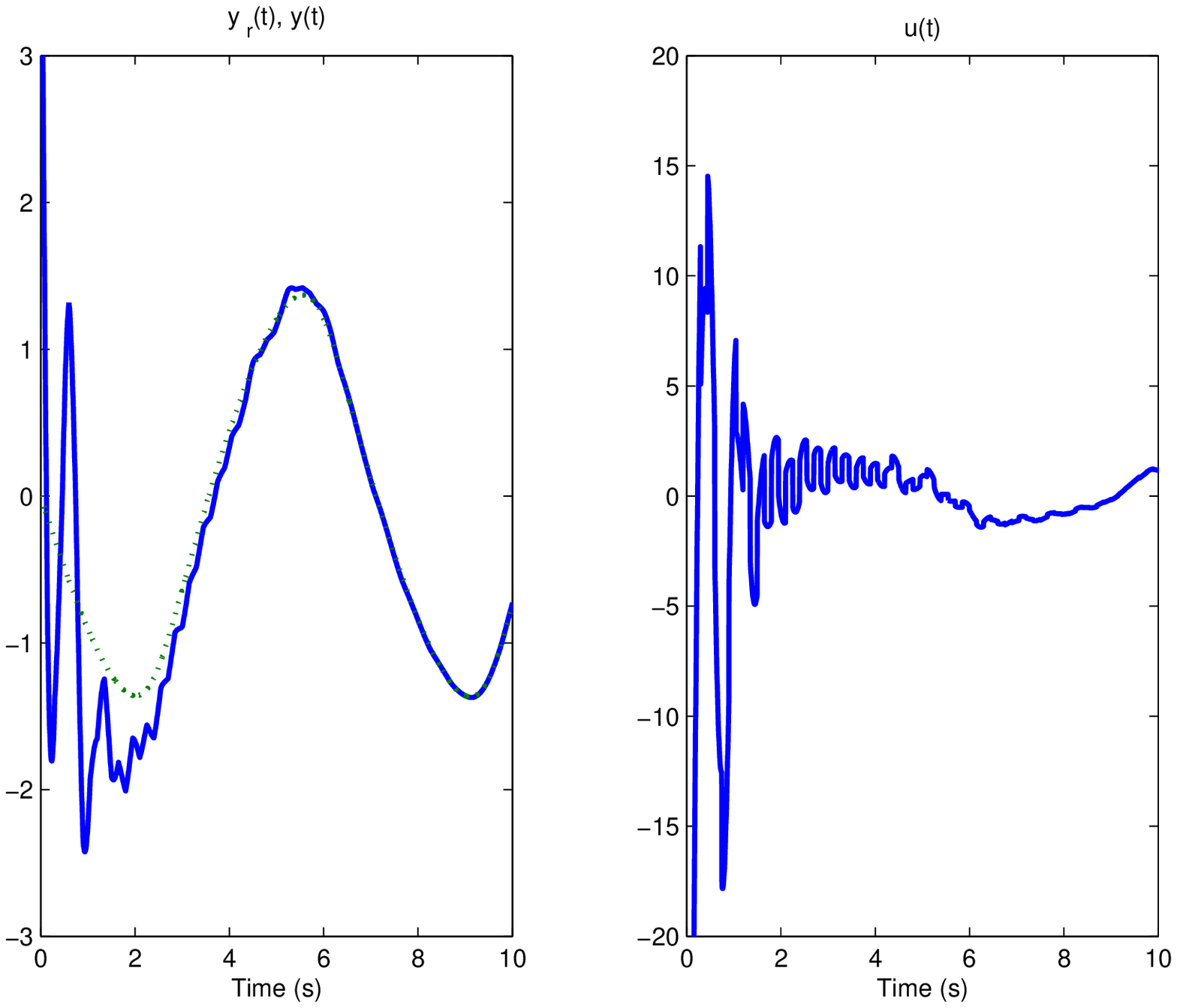}
\caption{\label{1cn3} First control scenario ($N=2$, $T=0.15$ s). Left: time
behavior of the reference trajectory $y_{\rm r}(t)$ (dotted line) and of the
controlled output $y(t)$ (solid line). Right: time behavior of the control
input $u(t)$.}
\end{figure}

\begin{figure}[htb]
\centering
\includegraphics[width=1\textwidth,height=0.5\textwidth]{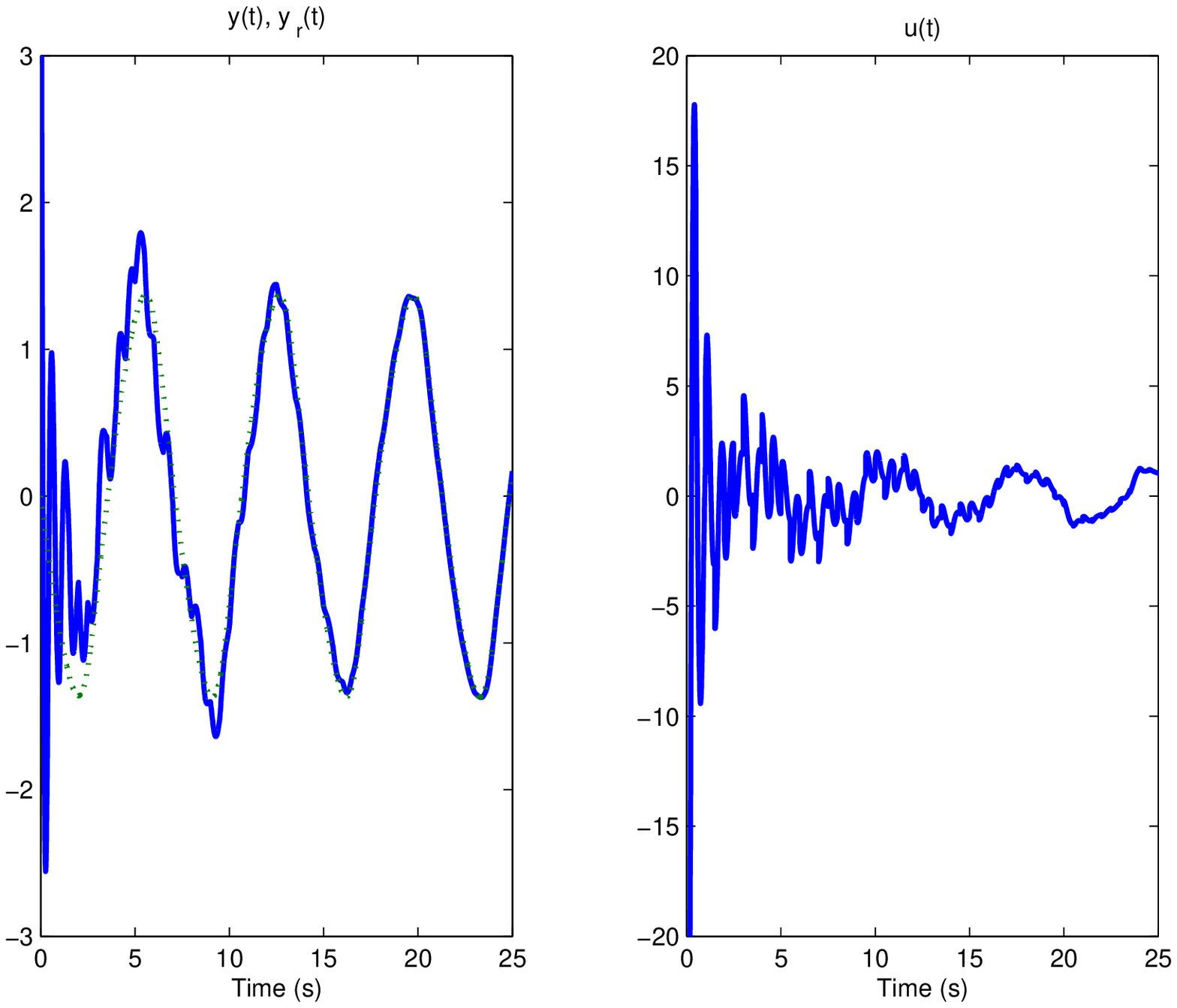}
\caption{\label{2cn3} Second control scenario ($N=4$, $T=0.5$ s). Left: time
behavior of the reference trajectory $y_{\rm r}(t)$ (dotted line) and of the
controlled output $y(t)$ (solid line). Right: time behavior of the control
input $u(t)$.}
\end{figure}

  We consider a Van der Pol oscillator with $\epsilon=1.5$ and
  $a=1$. The regulator (\ref{controllo1}) is tuned choosing
  $\kappa=3$, $G=\left (12 \,\, 36\right )\tr$ and $k=8$.
  We consider two different simulative scenarios which differ for the severity
  of the communication channel constraint. In the first case we
  suppose that the number of available bits is $N_b=2$
  yielding, according to (\ref{nbNrel}) and to the fact that $r=2$,  $N=2$. In
  this case, for a certain set of initial conditions, condition (\ref{TNconst})
  is fulfilled with $T=0.15$ s. In the second case the available number of
  bits is assumed $N_b=4$ from which (\ref{nbNrel}) and (\ref{TNconst})
  yield a bigger $N$ and $T$ respectively equal to $N=4$ and $T=0.5$ s. The
  simulation results, obtained assuming the exosystem (\ref{exoexa}) and the system
  (\ref{sysexa}) respectively at the initial conditions $w(0) = \left ( 1, 0 \right )$
  and $y(0)=5$, are shown in the figures \ref{1cn1}-\ref{1cn3} for the first
  scenario and figures \ref{2cn1}-\ref{2cn3} for the second one. In
  particular figure \ref{1cn1} (respectively \ref{2cn1}) shows the
  quantized variable $w_q$ transmitted from the encoder to the decoder and used to
  reset the respective dynamics according to the rule described in Section
  \ref{secEncDec}. Note that in the first scenario represented in figure \ref{1cn1}
  each of the two components of the vector $w_q$, taking value in the set
  $\{ -{1/2}, \,1/2\}$, can be transmitted using $1$ bit. On the other
  hand in the second scenario, represented in figure \ref{2cn1}, the transmission
  of each component of $w_q$, whose value are in the set $\{ -{3/2}, \,-1/2, \, 1/2, \, 3/2\}$
  requires $2$ bits.
  Figure \ref{1cn2} (respectively \ref{2cn2})  shows in the left-half side
  the error, as a  function of time, between the exosystem and the encoder
  (decoder) state and in the right-half side the phase portrait of the Van der Pol
  oscillator with overlapped the actual state trajectory of the encoder
  (decoder). Finally figure \ref{1cn3} (respectively \ref{2cn3}) plots the
  the tracking error $e(t) = y(t) - y_{\rm r}(t)$ on the left-half side and the
  control input $u(t)$ on the right-half side from which it is possible to see
  that in both the control scenarios the synchronization between (\ref{exoexa})
  and (\ref{sysexa}) is achieved.

\section{Conclusions}
We have discussed the problem of asymptotically tracking a reference signal
which is generated by a remote exosystem, and transmitted through a finite
bandwidth communication channel. Although only an estimate
of the actual tracking error is available to the regulator,
a suitable choice of the controller parameters
allows us to
achieve the control goal while fulfilling the constraint on the
bandwidth of the channel.

\end{document}